\documentclass[french,brochure]{bettyb}

\usepackage{url}

\usepackage{hyperref}

  \DeclareMathSymbol{A}{\mathalpha}{operators}{`A}%
  \DeclareMathSymbol{B}{\mathalpha}{operators}{`B}%
  \DeclareMathSymbol{C}{\mathalpha}{operators}{`C}%
  \DeclareMathSymbol{D}{\mathalpha}{operators}{`D}%
  \DeclareMathSymbol{E}{\mathalpha}{operators}{`E}%
  \DeclareMathSymbol{F}{\mathalpha}{operators}{`F}%
  \DeclareMathSymbol{G}{\mathalpha}{operators}{`G}%
  \DeclareMathSymbol{H}{\mathalpha}{operators}{`H}%
  \DeclareMathSymbol{I}{\mathalpha}{operators}{`I}%
  \DeclareMathSymbol{J}{\mathalpha}{operators}{`J}%
  \DeclareMathSymbol{K}{\mathalpha}{operators}{`K}%
  \DeclareMathSymbol{L}{\mathalpha}{operators}{`L}%
  \DeclareMathSymbol{M}{\mathalpha}{operators}{`M}%
  \DeclareMathSymbol{N}{\mathalpha}{operators}{`N}%
  \DeclareMathSymbol{O}{\mathalpha}{operators}{`O}%
  \DeclareMathSymbol{P}{\mathalpha}{operators}{`P}%
  \DeclareMathSymbol{Q}{\mathalpha}{operators}{`Q}%
  \DeclareMathSymbol{R}{\mathalpha}{operators}{`R}%
  \DeclareMathSymbol{S}{\mathalpha}{operators}{`S}%
  \DeclareMathSymbol{T}{\mathalpha}{operators}{`T}%
  \DeclareMathSymbol{U}{\mathalpha}{operators}{`U}%
  \DeclareMathSymbol{V}{\mathalpha}{operators}{`V}%
  \DeclareMathSymbol{W}{\mathalpha}{operators}{`W}%
  \DeclareMathSymbol{X}{\mathalpha}{operators}{`X}%
  \DeclareMathSymbol{Y}{\mathalpha}{operators}{`Y}%
  \DeclareMathSymbol{Z}{\mathalpha}{operators}{`Z}%

\iffalse
  \usepackage{fontspec}
  \setromanfont[Ligatures=TeX]{TeX Gyre Pagella}
  \usepackage{unicode-math}
\else
  \usepackage[utf8]{inputenc}
  \usepackage{newpxtext,newpxmath}
\fi

\usepackage[T1]{fontenc}
\usepackage[french]{babel}
\let\mathscr \mathcal

\author{Antoine Chambert-Loir}
\address{Univ. Paris Diderot, Sorbonne Paris Cité \\
 Institut de Mathématiques de Jussieu-Paris Rive Gauche, \\
 F-75013, Paris, France}
\email{Antoine.Chambert-Loir@math.univ-paris-diderot.fr}
\title{Le théorème de réduction stable de Deligne et Mumford}

\date{Janvier 2019}
\bbkannee{2\textsuperscript{e} ann\'ee, 2018--2019}
\bbknumero{7}

\begin{abstract}
L'espace de modules des courbes projectives lisses de genre~$g$ 
est une variété algébrique quasi-projective, mais non projective.
Pour comprendre sa géométrie, il est parfois crucial d'en considérer 
des compactifications.  En acceptant de paramétrer également des courbes 
(dites stables) aux singularités
contrôlées, Deligne et Mumford en ont construit une compactification projective.
Le caractère propre de cette compactification se traduit par le théorème de réduction stable qu'ils démontrent également. 
(Sa projectivité est un théorème ultérieur de Knudsen et Mumford.)
Cet texte, qui reprend l'exposé oral, est une introduction à ces objets.
\end{abstract}

\alttitle{The stable reduction theorem of Deligne and Mumford}
\begin{altabstract}
The moduli space of smooth projective curves of genus~$g$ is a
quasi-projective algebraic variety, but is not projective.
To understand its geometry, it may be crucial to consider compactifications
of this space. By allowing to parameterize as well curves with controlled
singularities (the so called stable curves), Deligne
and Mumford constructed a projective compactification.
The properness of this compactification translates into
the stable reduction theorem that they prove, its projectivity
is a later theorem of Knudsen and Mumford. This text is based
on the oral presentation and aims at introducing these objects.
\end{altabstract}

\def\C{\mathbf C}
\def\Z{\mathbf Z}
\def\P{\mathbf P}
\def\pgcd{\operatorname{pgcd}}
\def\Spec{\operatorname{Spec}}
\def\SL{\operatorname{SL}}
\def\PGL{\operatorname{PGL}}
\def\Card{\operatorname{Card}}

\begin{document}
\maketitle

\section{Surfaces de Riemann}

Les objets de base auxquels on s'intéresse ici sont
les surfaces de Riemann compactes (connexes, non vides).
Du point de vue \emph{topologique}, leur  classification est 
très simple, puisqu'elles sont distinguées par leur \emph{genre}
qui est un entier positif ou nul. Pour $g=0$, il s'agit de la sphère
$\mathbf S_2$; pour $g=1$, il s'agit d'un tore $(\mathbf S_1)^2$;
pour $g\geq 2$, on obtient des « tores à $g$ trous ».
Cependant, les surfaces de Riemann 
sont plus que leur espace topologique sous-jacent,
puisqu'elles sont munies d'une \emph{structure complexe} supplémentaire
et, en tant que telles, leur classification est plus délicate.
Le genre apparaît alors via le théorème de Riemann--Roch:
c'est la dimension de l'espace des formes différentielles holomorphes,
ou bien, avec un peu d'anachronisme, la dimension du premier groupe de cohomologie cohérente:
\[ g(X) = \dim ( H^0(X,\Omega^1_X)) = \dim (H^1(X,\mathcal O_X)). \]
Aussi, une forme différentielle méromorphe (non identiquement
nulle) possède $2g-2$ zéros comptés avec multiplicités,
les éventuels pôles étant comptés négativement:
\[ \deg( \Omega^1_X) = 2g-2. \]

Pour $g=0$, il n'y a qu'un modèle, la droite projective complexe $\mathbf P_1(\mathbf C)$, que l'on obtient en adjoignant au plan complexe $\C$ un point à l'infini, et la notion de fonction holomorphe au voisinage de l'infini étant donnée par le changement de variable $z'=1/z$.
Pour $g=1$, on parle de courbes elliptiques; elles s'obtiennent comme quotient du plan complexe~$\C$ par un réseau $\Lambda$, et deux tels quotients $\C/\Lambda$ et $\C/\Lambda'$ donnent lieu à la même courbe elliptique s'il existe
un nombre complexe non nul~$a$ tel que $\Lambda' = a\Lambda$. À homothétie
près, on peut ainsi choisir une base de~$\Lambda$ de la forme $(1,\tau)$, 
où $\tau$ est un élément du demi-plan de Poincaré~$\mathfrak h$,
et la classification des courbes elliptiques revient alors à l'étude du quotient
du demi-plan de Poincaré par l'action du groupe $\SL(2,\Z)$ agissant par homographies.

C'est le cas des surfaces de Riemann de genre~$g\geq 2$ qui fait l'objet
de cet exposé, et nous adopterons un point de vue algébrique:
d'après le théorème d'existence de Riemann, toute surface de Riemann 
compacte connexe peut être considérée comme une courbe projective lisse connexe
sur le corps des nombres complexes. De nombreux aspects de leur étude
s'étendent au cas des courbes sur un corps arbitraire, et même sur un anneau
arbitraire; le théorème de  Deligne et Mumford évoqué par le titre
de cet exposé en est un exemple remarquable.

Riemann avait déjà observé qu'une telle surface dépend de $3g-3$ 
paramètres complexes,
qu'il appelle \emph{modules}. 
Fricke et Teichmüller ont étudié plus précisément
l'ensemble des structures complexes sur une surface topologique orientée,
modulo isotopies ; l'\emph{espace de Teichmüller} qu'on obtient alors 
est une boule de l'espace~$\C^{3g-3}$.
Dans la théorie des déformations introduites par Kuranishi, l'entier
$3g-3$ apparaît comme la dimension de l'espace des formes différentielles
holomorphes de degré~$2$:
\[ 3g-3 = \dim(H^0(X,(\Omega^1_X)^{\otimes 2}). \]

Mais c'est à Mumford qu'on doit la construction de « l'espace des courbes
de genre $g\geq 2$ » : une variété algébrique complexe~$M_g$ 
munie d'une bijection de l'ensemble des classes d'isomorphie
de courbes projectives de genre~$g$  sur l'ensemble des points
complexes $M_g(\C)$ de cette variété, bijection qui est naturelle au sens suivant: pour toute famille $\mathcal C\to S$ de courbes projectives de genre~$g$,
l'application de~$S(\C)$ dans~$M_g(\C)$ donnée par $s\mapsto [\mathscr C_s]$
est régulière.
Cela demande de préciser la notion de \emph{famille}: il s'agit ici d'un 
morphisme propre et lisse dont toutes les fibres $\mathcal C_s$ en tout
point $s\in S(\C)$ sont des courbes de genre~$g$.
L'espace~$M_g$ est lisse, de dimension~$3g-3$,
quasi-projectif. 
Si $d$ est un entier~$\geq g+1$,
Severi a démontré que
toute courbe de genre~$g$
peut être vue comme un revêtement de degré~$d$ de la sphère de Riemann, 
dont la ramification est le plus simple possible
(au plus un point double dans chaque fibre), puis,
utilisant les travaux de Hurwitz et Clebsch, en a déduit que \emph{$M_g$
est connexe.}

En fait, Mumford a construit un schéma~$M_g$ sur~$\Z$,
qui paramètre les courbes algébriques de genre~$g$ sur tout corps
algébriquement clos. La connexité de ses fibres
est le résultat principal de l'article de Deligne et Mumford.
Fulton avait pu étendre l'approche de Severi à la caractéristique assez grande.

Pour la construction  de Mumford,  deux propriétés géométriques d'une
courbe~$X$ de genre~$g\geq 2$ sont fondamentales:
\begin{enumerate}
\item Comme $3(2g-2)\geq 2g+1$, 
le fibré tricanonique $(\Omega^1_X)^{\otimes 3}$ est très ample;
\item Le groupe d'automorphismes de~$X$ est fini
(précisément, on a l'inégalité de Hurwitz $\Card(\operatorname{Aut}(X))\leq 42(2g-2)$.
\end{enumerate}

Comme 
\[ \dim(H^0(X, (\Omega^1_X)^{\otimes 3}))=3(2g-2)+1-g=5g-5, \]
la première propriété entraîne que toute courbe de genre~$g\geq 2$
peut être plongée dans l'espace projectif $\P_{5g-6}$
et permet d'appliquer 
la théorie des schémas de Hilbert, due à Grothendieck, 
qui construit ainsi l'ensemble des sous-schémas fermés de polynôme
de Hilbert donné d'un espace projectif.
En effet, le polynôme de Hilbert de l'image~$\iota(X)$ est déterminé:
pour tout entier~$n\geq 1$, on a
\begin{align*}
 P_{\iota(X)} (n)  & = \dim (H^0(\iota(X), \mathscr O(n))) = \dim(H^0(X,(\Omega^1_X)^{\otimes 3n})) \\
& =3n(2g-2)+1-g=(6n-1)(g-1). \end{align*}
On peut alors démontrer que $\iota(X)$ est définie par des équations
homogènes dont le degré est contrôlé
% \footnote{Indiquer ce qu'on trouve}
et représenter $\iota(X)$ par le sous-espace de ses équations de tel degré
parmi l'espace des polynômes homogènes de ce même degré.
Le schéma de Hilbert~$H_g$ des
courbes tricanoniquement plongées apparaît ainsi 
comme sous-schéma d'une variété grassmannienne.

La seconde étape de la construction consiste à « oublier » le choix
d'un plongement, c'est-à-dire à passer au quotient l'espace~$H_g$
par l'action naturelle du groupe $\PGL(5g-5,\C)$.
C'est là qu'intervient la \emph{théorie géométrique des invariants}
développée par Mumford.
Il démontre en effet que l'algèbre graduée des « polynômes » sur~$H_g$
qui sont invariants sous l'action du groupe $\PGL(5g-5,\C)$
est de type fini, et que des points de~$H_g$ qui ne sont
pas dans la même orbite sont distingués
par ces polynômes invariants.
(C'est là qu'intervient la finitude du groupe d'automorphismes.)
L'espace de modules~$M_g$ n'est alors autre qu'un ouvert du
spectre projectif de cette algèbre graudée.

Par construction, l'espace~$M_g$ ainsi construit est quasi-projectif.
Pour en comprendre les propriétés géométriques de~$M_g$,
il est nécessaire de le compactifier d'une façon naturelle.
La façon qu'ont proposée Deligne et Mumford consiste à
lui adjoindre des courbes qui soient assez singulières
(pour ajouter assez de points), mais pas trop (pour que la compactification
obtenue reste raisonnable); c'est la notion de \emph{courbe stable}.

\section{Courbes stables}

Deligne et Mumford appellent \emph{courbe stable} sur un corps~$k$
algébriquement clos un schéma~$C$ propre, connexe, 
réduit, de dimension~$1$ vérifiant les deux propriétés suivantes:
\begin{enumerate}
\item Ses singularités sont des points doubles ordinaires;
\item Si $E$ est une composante irréductible lisse de~$C$
de genre~$0$, alors $E$ rencontre les autres composantes en au moins
trois points.
\end{enumerate}
Il y a plusieurs façons de donner un sens précis à la première 
propriété.  On peut la formuler pour la « topologie étale » :
$C$ est alors localement croisement de deux axes de coordonnées.
On peut aussi décrire les anneaux locaux complétés :
pour tout point singulier $c\in C$, il existe un isomorphisme
de $k[[x,y]]/(xy-a)$ sur $\widehat{\mathcal O_{C,c}}$.
On appelle alors \emph{genre} d'une telle courbe stable~$C$ la dimension
du premier groupe de cohomologie cohérente:
\[ g(C) = \dim_k (H^1(C,\mathcal O_C)). \]

On observe que cette notion est stable par extension de
corps algébriquement clos.

Les singularités d'une courbe stable~$C$ sont les plus simples possibles ; 
elles garantissent l'existence d'un \emph{faisceau dualisant}~$\omega_C$ 
qui est localement libre de rang~$1$.
Sa description est également très simple.
Notons $p\colon C'\to C$ la normalisation de~$C$;
si $z_1,\dots,z_n$ sont les points singuliers de~$C$,
notons $x_i,y_i$ les deux antécédents de~$z_i$ par~$p$. Alors,
le faisceau $p^{-1}(\omega_C)$ s'identifie au sous-faisceau
de $\Omega^1_{C'}(\sum x_i + \sum y_i)$, faisceau des formes différentielles
sur~$C'$ à pôles au plus simples en les~$x_i$ et~$y_i$,
formé des formes différentielles~$\omega$
telles que, pour tout~$i$, la somme des résidus de~$\omega$
en~$x_i$ et~$y_i$ soit nulle.

Grâce à la seconde propriété de la définition d'une courbe stable,
on démontre que $\omega_C^{\otimes n}$ est très ample pour $n\geq 3$
et plonge~$C$ dans~$\P_{5g-6}$; son image est un sous-schéma fermé
de polynôme de Hilbert $(6n-1)(g-1)$.
De plus, le groupe d'automorphismes d'une courbe stable~$C$
est fini.

Deligne et Mumford construisent alors un espace de modules $\overline M_g$
des courbes stables de genre~$g$ sur~$\Z$.
C'est un schéma de type fini sur~$\Z$ muni, pour tout corps algébriquement
clos~$k$, d'une bijection de l'ensemble des classes d'isomorphie
de courbe stable de genre~$g$ sur l'ensemble~$\overline M_g(k)$,
de sorte que pour toute famille $\mathscr C\to S$ de courbes
stables de genre~$g$, l'application $s\mapsto [\mathscr C_s]$
de~$S$ dans~$\overline M_g$ soit régulière.
Là aussi, il faut préciser  la notion de \emph{famille}: il s'agit
d'un morphisme $\mathscr C\to S$ propre et plat 
dont les fibres géométriques soient des courbes stables de genre~$g$.
Plus explicitement : 
pour tout point~$s$ de~$S$, le schéma $\mathscr C_{s}$ sur~$\overline{\kappa(s)}$ soit une courbe stable, où $\kappa(s)$ est une clôture algébrique
du corps résiduel~$\kappa(s)$ de~$s$. 
Toute courbe lisse de genre~$g$ est évidemment stable, de sorte que
l'espace~$M_g$ apparaît comme une partie de~$\overline M_g$,
ouverte et dense.

C'est une compactification de l'espace des modules des courbes.
En fait, Deligne et Mumford construisent un objet plus général
et plus précis que $\overline M_g$ dans lequel on garde une trace
des automorphismes des courbes stables : c'est le \emph{champ des modules}
$\overline{\mathcal M}_g$ des courbes stables de genre~$g$.
Il y a (dans la catégorie des champs algébriques) un morphisme canonique
$\overline{\mathcal M}_g\to \overline{M}_g$
qui fait de l'espace de modules des courbes stables
l'espace grossier de ce champ.  

\emph{Le champ $\overline{\mathcal M}_g$ est propre et lisse sur~$\Z$,
de dimension relative $3g-3$,
à fibres géométriques connexes} (Deligne, Mumford). 
\emph{L'espace de modules grossier
$\overline{M}_g$ est projectif} (Gieseker, Knudsen, Mumford).
On peut également décrire précisément le complémentaire
 $\overline{\mathcal M}_g\setminus\mathcal M_g$ 
du champ des courbes lisses dans celui des courbes stables,
ou bien $\overline{M}_g\setminus M_g$ de l'espace grossier
des courbes lisses dans celui des courbes stables.

Une fois la construction de ce champ acquise, sa lissité découle
de la théorie des déformations des courbes stables. Sa propreté
sera discutée au paragraphe suivant,
et la connexité de ses fibres en découle, car on sait déjà qu'il 
en est ainsi de sa fibre au-dessus du point générique de~$\Spec(\Z)$,
par la théorie classique.

\section{Le théorème de réduction stable}

On l'a déjà annoncé au paragraphe précédent, le champ~$\overline{\mathcal M}_g$
des courbes stables de genre~$g$ est \emph{propre}.
Le théorème de réduction stable de Deligne et Mumford
est  la traduction géométrique du critère valuatif de propreté
pour ce champ.

La partie « existence » de ce critère est le théorème suivant: \emph{Soit $R$ un anneau de valuation discrète et soit $K$ son corps des fractions, soit $C$ une courbe projective lisse connexe de genre~$g\geq 2$ sur~$K$. Il existe une extension finie~$K'$ de~$K$
et, notant~$R'$ la clôture intégrale de~$R$ dans~$K'$,
une famille stable $\mathcal C'\to \Spec(R')$ de fibre
générique~$\mathcal C'_{K'}=C_{K'}$.}

La partie « unicité » affirme, quant à elle, que \emph{si $\mathcal C$
et $\mathcal C'$ sont des $S$-courbes stables, le faisceau
$\operatorname{Isom}_S(\mathcal C,\mathcal C')$ de leurs
isomorphismes est fini et non ramifié sur~$S$.}
En particulier, dans les conditions du théorème de réduction
stable, la famille $\mathcal C'$ est déterminée par sa fibre générique,
c'est-à-dire par la courbe~$C$.

L'importance du théorème de Deligne et Mumford, et sa difficulté,
réside pour une grande part dans sa généralité ; en particulier,
le corps résiduel de l'anneau de valuation discrète~$R$
peut être de caractéristique arbitraire, et c'est la clef pour
en déduire l'irréductibilité en toute caractéristique
de l'espace des modules des courbes de genre~$g$.

Nous allons cependant voir qu'en caractéristique zéro, 
ce théorème possède une démonstration assez simple.

\medskip

Tout d'abord, nous allons raisonner comme si l'anneau de valuation
discrète~$R$ était l'anneau des fonctions holomorphes
sur un disque~$\Delta=D(0;1)$ du plan complexe, la courbe~$C$ sur~$K$
étant donnée par une surface~$\mathscr S^*$ fibrée
sur le disque épointé~$\Delta^*=\Delta\setminus\{0\}$ dont les fibres
sont des courbes lisses de genre~$g$.

Une première étape consiste à « prolonger » cette surface
au-dessus du point manquant~$0$: on part ainsi d'une surface~$\mathscr S$
munie d'un morphisme propre et plat $p\colon \mathscr S\to\Delta$
dont les fibres $\mathscr S_t$, pour $t\in\Delta^*$,
sont des courbes lisses de genre~$g$, mais dont la fibre~$\mathscr S_0$ 
est éventuellement singulière.
En utilisant le théorème de résolution des singularités des surfaces,
on suppose que $\mathscr S$ est lisse; en utilisant également
le théorème de résolution plongée des courbes, on se ramène
au cas où les composantes irréductibles de~$\mathscr S_0$ 
sont lisses et se coupent transversalement.
Le diviseur de la fonction holomorphe~$t$ sur~$\mathscr S$
est de la forme $\mathscr S_0=\sum n_i E_i$: la famille $(E_i)$ est la famille
des composantes irréductibles de~$\mathscr S_0$ et pour tout~$i$,
$n_i$ est la multiplicité de~$E_i$ dans~$\mathscr S_0$.

Cette première étape pourrait fonctionner en caractéristique positive,
ou en caractéristique mixte,
car on y dispose des théorèmes de résolution des singularités.
Une seconde étape vise à \emph{éliminer les multiplicités}
et la méthode que nous allons utiliser ne fonctionnerait pas
en caractéristique positive.

Notons $n$ le ppcm des~$n_i$ et considérons le revêtement
ramifié $\Delta'=\Delta\to \Delta$
donné par $t\mapsto t^{n}$. Soit $\mathscr S'$ la surface
déduite de~$\mathscr S$ par ce changement de base;
explicitement, c'est l'ensemble des couples $s'=(s,t)\in\mathscr S\times\Delta'$
tels que $p(s)=t^n$, considéré comme famille au-dessus de~$\Delta'$
par la seconde projection $(s,t)\mapsto t$.
Autour d'un point~$s'$ de la forme $(s,t)$ avec $t\neq 0$,
la situation est aussi belle qu'auparavant, car le revêtement
$\Delta'\to \Delta$ n'est ramifié qu'au-dessus de l'origine.
Pour analyser la structure locale de la surface~$\mathscr S'$
autour d'un point~$s'$ de la forme $(s,0)$, il y a deux
cas, suivant que $s$ appartient à une seule, ou à deux composantes
de~$\mathscr S_0$.

a) Supposons d'abord que $s$ n'appartient qu'à une seule composante
de~$\mathscr S_0$, de multiplicité~$a$.
Alors, on peut décrire~$\mathscr S$ au voisinage de~$s$
par deux paramètres $(u,v)$ et la projection vers~$\Delta$
étant donnée par $(u,v)\mapsto u^a$. 
Au voisinage de~$s'$, on décrit alors~$\mathscr S'$
comme l'ensemble des triplets $(u,v,t)$ tels que $u^a=t^n$ et
cette équation  se factorise en $\prod_{\zeta^a=1} (u-\zeta t^m))=0$, 
où $m=n/a$. On obtient une réunion de $a$~surfaces lisses
(paramétrées par $u=\zeta t^m$, pour $\zeta^a=1$), de projection
vers~$\Delta$ données par $(u,v,t)\mapsto t$; leur fibre en~$0$
est d'équation $u=0$: il n'y a plus de multiplicité.

b) Supposons maintenant que $s$ appartient à deux composantes
de~$\mathscr S_0$, de multiplités~$a$ et~$b$. 
On peut encore décrire~$\mathscr S$ au voisinage de~$s$ par deux
paramètre~$(u,v)$ mais la projection vers~$\Delta$ est maintenant
donnée par $(u,v)\mapsto u^a v^b$. Posons $d=\pgcd(a,b)$,
soit $\alpha$ et $\beta$ des entiers (premiers entre eux)
tels que $a=d\alpha$ et $b=d\beta$.
Soit enfin  $m=n/d \alpha\beta$; c'est un entier.
Le changement de base 
fournit la surface de~$\C^3$ d'équation $u^{d\alpha} v^{d\beta}=t^{md\alpha\beta}$;
elle se décompose comme la réunion des $d$~surfaces
d'équations $u^{\alpha} v^{\beta}= \zeta t^{m\alpha\beta}$, 
où $\zeta$ parcourt
les racines $d$-ièmes de l'unité. On est ainsi essentiellement
ramené au cas où $a$ et $b$ sont premiers entre eux.

L'hypersurface de~$\C^3$ d'équation $u^av^b=t^{abm}$ est singulière
en l'origine. Démontrons que sa normalisation
est l'hypersurface d'équation $xy=t^m$, par
l'application $(x,y,t)\mapsto (x^b,y^a,t)$.
Cette application est en effet génériquement bijective:
soit $(x_1,y_1,t_1)$ et $(x_2,y_2,t_2)$ des points de~$\C^3$
qui vérifient $x_1y_1=t_1^m$, $x_2y_2=t_2^m$ et $x_1^b=x_2^b$, $y_1^b=y_2^a$
et $t_1=t_2$;
posons $\zeta=x_2/x_1$ et $\xi=y_2/y_1$; on a $\zeta^b=\xi^a=\zeta\xi=1$,
donc $\zeta=\xi=1$ puisque $a$ et $b$ sont premiers entre eux.
L'hypersurface d'équation $xy=t^m$ est singulière en l'origine,
avec une singularité normale, de type~$A_{m-1}$. Sa résolution explicite
est bien connue, elle insère « dans le point singulier »
une chaîne de $(m-1)$~droites projectives, toutes de multiplicité~$1$.

À l'issue de cette deuxième étape, 
nous avons, au prix d'un revêtement ramifié du disque,
prolongé la surface~$\mathcal S^*$ en une surface lisse~$\mathcal S'$ 
munie d'un morphisme vers~$\Delta'$ dont la fibre centrale
est réunion de composantes lisses, de multiplicités~1, qui se coupent
transversalement. Pour obtenir une famille stable, il reste à assurer
la seconde condition de la définition: que les composantes lisses
de genre~$0$ coupent les autres en au moins trois points,
ce qui est l'objet de la troisième étape.

Soit $(E_i)$ la famille des composantes irréductibles de~$\mathcal S'_0$;
on peut écrire une égalité de diviseurs $\mathcal S'_0=\sum E_i$.
La théorie de l'intersection sur la surface~$\mathcal S'$
permet d'écrire, si $E$ est une composante,
\[ - E^2 = E \cdot (\mathcal S'_0-E) \]
est le nombre de points d'intersection de~$E$ avec les autres composantes.
Si $E$ contredit la définition d'une famille stable, on a donc
$E^2\in\{0,-1,-2\}$. Le cas $E^2=0$ est impossible;
il signifierait que $E$ est isolée
dans sa fibre, ce qui contredit le \emph{Main theorem} de Zariski.
De même, ces courbes exceptionnelles se regroupent en chaînes
et la théorie des surfaces permet de les contracter.
On démontre que la surface~$\mathcal S''$ qui en résulte
fournit la famille stable cherchée.

\section{Et après ?}

La combinaison des travaux de nombreux mathématiciens
effectués dans les vingt dernières années,
notamment ceux de Kollár, Hacon, McKernan et Xu, 
a permis de généraliser ce théorème de réduction stable aux familles
de variétés de dimension strictement supérieure à~$1$,
en caractéristique nulle.

Une première étape, assez peu différente du cas des courbes,
consiste à construire une famille dont l'espace total est lisse
et dont la fibre centrale est un diviseur à croisements normaux.
On y parvient grâce au théorème de résolution des singularités
d'Hironaka.

Une seconde étape permet de supposer que les composantes
de la fibre centrale sont de multiplicité~$1$. 
Il s'agit là du théorème de réduction \emph{semi-stable}
de Kempf, Knudsen, Mumford et Saint-Donat.

La dernière étape est beaucoup plus difficile et demande déjà
d'avoir précisé la définition de \emph{variété stable} 
en dimension supérieure: 
ce seront  des variétés aux singularités \emph{semi log-canoniques}
pour lesquelles le faisceau canonique est ample.
Son aboutissement dépend des progrès récents
dans la théorie du programme des modèles minimaux (\textsc{mmp}).

Je renvoie à l'exposé d'Olivier Benoist au Séminaire Bourbaki
pour une introduction à ces développements.

% initié par Mori, notamment le théorème, dû à Hacon et McKernan,
% de finitude de l'anneau canonique.

\bibliographystyle{smfplain}
\bibliography{aclab,acl,semis}
\nocite{deligne-m69,kempf-knudsen-mumford-saintdonat1973,knudsen1983b,knudsen1983c,knudsen-mumford1976}
\nocite{arbarello-cornalba-griffiths2011,mumford-f-k94,grothendieck1961,benoist2019}
\end{document}